# TR-Based Antenna Design with Forward FD: the Effects of Step Size on the Optimization Performance

Adrian Bekasiewicz[1][0000-0003-0244-541X], Slawomir Koziel[1,2][0000-0002-9063-2647],
Tom Dhaene[3][0000-0003-2899-4636], and Marcin Narloch[1][0000-0002-7640-2941]

[1] Faculty of Electronics, Telecommunications and Informatics, Gdansk University of Technology, Narutowicza 11/12, 80-233 Gdansk, Poland
[2] Department of Engineering, Reykjavik University, Menntavegur 1, 102 Reykjavík, Iceland
[3] Department of Information Technology (INTEC), IDLab, Ghent University-imec, iGent, Technologiepark-Zwijnaarde 126, 9052 Ghent, Belgium
`bekasiewicz@ru.is`

**Abstract.** Numerical methods are important tools for design of modern antennas. Trust-region (TR) methods coupled with data-efficient surrogates based on finite differentiation (FD) represent a popular class of antenna design algorithms. However, TR performance is subject to FD setup, which is normally determined *a priori* based on rules-of-thumb. In this work, the effect of FD perturbations on the performance of TR-based design is evaluated on a case study basis concerning a total of 80 optimizations of a planar antenna structure. The obtained results demonstrate that, for the considered radiator, the performance of the final designs obtained using different FD setups may vary by as much as 18 dB (and by over 4 dB on average). At the same time, the *a priori* perturbations in a range between 1.5% and 3% (w.r.t. the initial design) seem to be suitable for maintaining (relatively) consistent and high-quality results.

## 1 Introduction

Optimization methods are indispensable tools for the design of contemporary antennas. The main bottleneck of conventional algorithms is that they require dozens, or even hundreds of model evaluations to converge [1], [2]. At the same time, the simulation cost of antenna electromagnetic (EM) models is high—especially for modern multi-parameter topologies—which challenges the concept of their direct optimization, especially if the designer does not have access to high-performance computing clusters and a large number of software licenses [3], [4]. From this perspective, availability of data-efficient methods is of high importance for reliable, yet low-cost antenna design.

Data-efficient design of microwave and antenna circuits can be performed using surrogate-assisted methods (SAM). The goal of SAM is to embed the antenna design into a meta-loop that involves optimization of a cheap auxiliary model (with limited accuracy) that is iteratively updated/re-constructed using only a handful of accurate (yet numerically expensive) data samples obtained from EM simulations [5], [6]. Trust-region (TR) methods belong to a popular class of SAM algorithms. They perform exploitation of the search space through optimization of the local, numerically cheap model that predicts the new design candidates. To alleviate the high cost of EM simulations, the model is often in the form of a first-order Taylor expansion constructed from the structure response and its derivatives [1], [7]. Although the latter can be obtained using analytical derivatives—at a low computational overhead (w.r.t. zero-order simulation)—adjoint-capable simulations are available in only a handful of commercially available EM solvers [8], [9]. Alternatively, the sensitivity data can be approximated using finite-differences (FD), i.e., through evaluation of the structure response at the design of interest and series of small perturbations around it. In antenna engineering, linear models are often constructed using forward FDs as which require only one EM simulation per design parameter [2], [8].

Regardless of proved usefulness, the performance of FD-based TR optimization is subject to appropriate setup (perturbation-wise). Although algorithms dedicated to automatically determine suitable FD steps have been reported in the literature, they are prohibitively expensive when applied to EM-driven design problems [10], [11]. Instead, the perturbations are predominantly selected *a priori* based on rules-of-thumb. Popular methods include determination of FD steps as a fraction of the initial design, square-root of machine precision, or as a based on experience-driven manual tuning of structure [7], [10]. Despite explicit relation between the FD setup and TR optimization performance, the problem pertinent to determination of appropriate perturbations for EM-driven optimization of antennas remains open.

In this work, the effects of *a priori* selected FD steps on the performance of TR-based gradient optimization have been evaluated on a case study basis concerning design of a planar quasi-patch antenna. A total of 80 numerical experiments spanning across 10 different designs and 8 FD-setups each has been performed. The results demonstrate that, for the considered antenna, the FD-induced performance discrepancy between the optimized designs is substantial and can vary by as much as 18 dB. At the same time, the average (i.e., calculated over all of the designs) best-to-worst-setup difference amounts to over 4 dB. The obtained data not only demonstrate that appropriate FD is crucial for ensuring high performance (and consistency of optimization) but also suggests that perturbations in a range between 1.5% to 3% of the initial design are suitable for mitigating the effects of numerical noise while supporting exploitation of the search space using the linear approximation models.

## 2   Optimization Algorithm

### 2.1   Problem Formulation

Let $R(x)$ be the EM-simulation response of the structure obtained for the vector of input parameters $x$. The optimization problem is given as:

$$x^* = \arg\min_{x \in X} U(R(x)) \qquad (1)$$

Here, the $x^* \in X$ represents the optimized design to be found within the feasible region of the search space $X$, which is defined by the lower/upper bounds $l/u$; $U$ is a scalar objective function. Due to high evaluation cost, direct optimization of $R$ is impractical. Instead, the task (1) can be embedded into a surrogate-assisted design framework to enable cost-efficient identification of $x^*$.

## 2.2   TR Optimization Framework

The goal of TR-based optimization is to generate a series of approximations, $i = 0, 1, 2, \ldots$, to the problem (1), as [6]:

$$x^{(i+1)} = \arg \min_{\delta^{(i)} \geq \|x - x^{(i)}\|} U\left(G^{(i)}(x)\right) \qquad (2)$$

where $G^{(i)} = R(x^{(i)}) + J(x^{(i)})(x - x^{(i)})$ is a linear approximation model obtained at the $i$th step of the optimization process and $J$ represents its FD-based Jacobian [10]:

$$J(x^{(i)}) = \begin{bmatrix} \frac{1}{p_1}\left(R(x^{(i)} + p_1) - R(x^{(i)})\right) \\ \vdots \\ \frac{1}{p_D}\left(R(x^{(i)} + p_D) - R(x^{(i)})\right) \end{bmatrix}^T \qquad (3)$$

The problem (2) is solved by a gradient-based algorithm executed on the linear approximation models constructed from (3). It should be noted that the variable $p_d$ ($d = 1, \ldots, D$) represents the perturbation w.r.t. $d$th dimension of the design $x^{(i)}$, whereas $p_d = [0 \ldots p_d \ldots 0]^T$. The trust-region radius $\delta^{(i)}$ is adjusted based on the gain coefficient $\rho = [U(R(x^{(i+1)})) - U(R(x^{(i)}))]/[U(G^{(i)}(x^{(i+1)})) - U(G^{(i)}(x^{(i)}))]$ which expresses the ratio between the expected and the obtained objective function change [6]. The factor is also used to accept ($\rho > 0$), or reject ($\rho < 0$) the candidate designs generated by (2) as new approximations of the final solution. The initial radius is set to $\delta^{(0)} = 1$. For the remaining iterations it is decreased as $\delta^{(i+1)} = \alpha_1\|x^{(i+1)} - x^{(i)}\|$ for $\rho < 0.05$ (poor performance of the local model), or increased as $\delta^{(i+1)} = \max(\alpha_2\|x^{(i+1)} - x^{(i)}\|, \delta^{(i)})$ for $\rho > 0.9$ (acceptable improvement of the response). The scaling factors are set to $\alpha_1 = 0.25$, and $\alpha_2 = 2.5$ [6]. The algorithm is terminated when $\delta^{(i+1)} < \varepsilon$ or $\|y^{(i+1)} - y^{(i)}\| < \varepsilon$ (here, $\varepsilon = 10^{-2}$). For more detailed discussion of the method, see [2], [6], [8].

## 2.3   Perturbations in FD-Based Jacobians

On the conceptual level, the problem pertinent to selection of appropriate perturbations (also referred to as a step-size dilemma; cf. Fig. 1(a)) involves balancing their size so as to minimize both the truncation (too large step) and round-off (too small step) errors [10]. Another challenge is that EM simulation models are inherently noisy which stems from inconsistency of their discretization for closely-located designs (as the ones required for construction of the linear models). As a consequence, for too small step the linear model might approximate noise, while the use of too large per-turbations might result in poor representation of the objective function changes (see Fig. 1(b)) [12], [13]. Furthermore, the problem of FD-setup is both circuit- and design-dependent. Consequently, a series of model evaluations would be required each candidate solution and w.r.t. each variable to find optimal steps. Clearly, it is not practical when the EM-driven antenna development is considered [10], [11], [14].

The FD setup for antenna design problems is normally performed through *a priori* specification of perturbations that are (hopefully) small enough to capture local changes of the objective function but also large enough to mitigate the effects of numerical noise on the identified descent direction (see Fig. 1(b)). The popular rule-of-thumb methods include selection of steps that either represent a fraction of the initial design, or correspond to the machine precision (single for the EM-driven problems) [10]. Alternatively FD-steps can be obtained based on the experience-driven tuning that involves visual inspection of the objective function changes [7]. Notwithstanding, the question concerning practical usefulness of the mentioned concepts for EM-based problems, as well as their effects on TR algorithm performance remains open.

## 3    Antenna Structure

Figure 2 illustrates the quasi-patch antenna considered for experiments [9]. The structure is dedicated to ISM-band (industrial, scientific, medicine) applications. It is designed on a dielectric substrate with permittivity/thickness of 3.5/0.762 mm and comprises a bi-component radiator in the form of a deformed patch loaded by a monopole strip and fed through a 50 Ohm microstrip line. Impedance matching of the patch and the feedline is maintained using two insets (notches). The ground plane length below the radiating component is parameterized so as to ensure its flexibility in terms of operational bandwidth (w.r.t. conventional structures). The component is implemented in CST Microwave Studio, and evaluated using its time-domain solver [15]. On average, the EM model is discretized using 150,000 tetrahedral mesh cells, whereas its evaluation cost amounts to 120s. The model also implements the coaxial connector that is normally used to mate the antenna with the test equipment [16].

The vector of design parameters is $x = [L\ l_2\ W\ w_2\ l_0\ o_0]^T$. The dependent variables are $o = 0.22L$, and $l_s = 0.1L$, whereas dimensions $l_1 = 1.5$, $w_1 = 2.5$, $w_s = 0.5$, and $w_0 = 1.7$ remain constant thorough the optimization process (all in mm). The feasible region of the design space $X$ is defined by the following lower and upper bounds: $l = [10\ 5\ 3.5\ 0.2\ 3\ 2]^T$ and $u = [25\ 25\ 10\ 3.2\ 15\ 10]^T$. The objective function $U(x) = U(R(x))$ is given as:

$$U(x) = \max\{R(x)\}_{f_L \leq f \leq f_H} \qquad (4)$$

where $f_L = 5$ GHz and $f_H = 6$ GHz represent the corner frequencies for the antenna operational bandwidth, whereas $R(x) = R(x, f)$ represents its reflection response and $f \in f$ is the frequency sweep. Note that the function (4) defines a simple min-max problem, where the goal is to minimize the maximum value of the antenna response within the frequency range of interest.

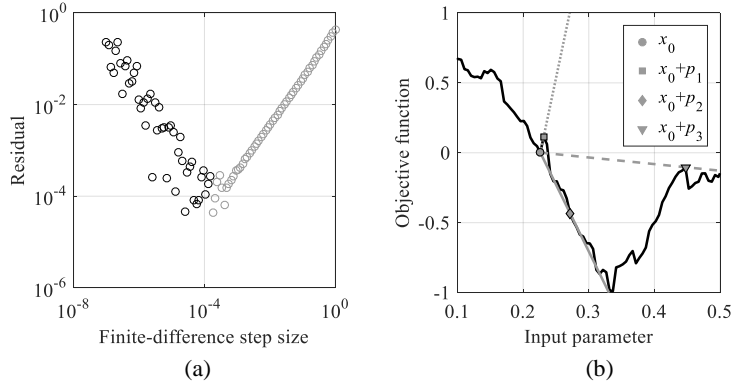

**Fig. 1.** FD steps: (a) a residual of analytical and FD-based derivatives for a function $\sin(\pi/4)$ with highlight on the round-off (black) and truncation errors (gray), (b) the quality of linear models for a noisy function when FD steps are too small ($\cdots$), too large (– –), and correct (—).

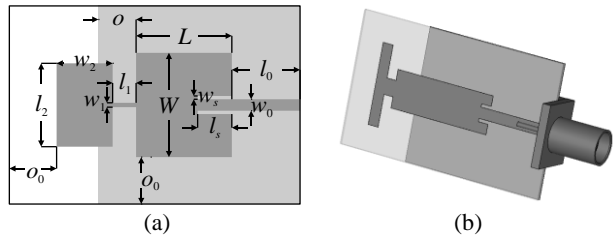

**Fig. 2.** Quasi-patch: (a) geometry with highlight on design variables and (b) visualization of the structure with coaxial connector [9]. Note that light gray represents the ground plane.

## 4   Numerical Results

The considered numerical experiments included a series of TR-based optimizations of the antenna for ten—randomly selected within the *l*/*u* bounds—designs (for details, see Table 1), each using eight different FD setups. These included specification of steps as: a manually defined fraction of the initial design (*p* from 0.5% to 3% with 0.5% step), or based on a single-point machine precision (double is not supported by the EM solver; $p = \varepsilon^{0.5} \approx 0.0032$; $\varepsilon = 10^{-7}$), as well as through manual tuning with $\boldsymbol{p} = [p_1 \ldots p_D] = 0.001 \cdot [3\ 3\ 3\ 7\ 8\ 6]^T$ (here, $D = 6$; cf. Sections 2 and 3). Note that the computational cost of experience-driven tuning (pertinent to model evaluations necessary for visual inspection of function changes) corresponds to 10 EM simulations.

The results gathered in Table 2 indicate that the FD perturbations setup substantially affects the TR convergence and the performance of the obtained designs. For the considered problems and tests, the maximum discrepancy between the designs obtained using different settings amounts to 18 dB (per specific design – $\boldsymbol{x}_7$) and 4.2 dB (on average for all of the designs). Furthermore, results suggest that, although the best responses (collectively) are obtained using 3% steps, for design $\boldsymbol{x}_9$ the performance of optimized responses was almost 4 dB lower compared to the case with 2.5% perturba-tions. The result suggests that the local model slightly overshoot the optimum. On the other hand, the calculated standard deviations $\sigma$ indicate that 3% step offers the most consistent responses among all experiments (the obtained value would be even lower if not for the poor-quality design $x_{10}$). As it comes to the test cases that involved custom and precision-based perturbations, their performance is comparable to 0.5% and 1% steps (which is understandable given their comparable ranges). Worsened results for small perturbations suggest that the driving factor behind the premature convergence of TR is the numerical noise (which surpasses round-off errors). From this perspective, the use of steps in the range from 1.5% to 3% of $x$ seem to counteract the noise while (mostly) ensuring appropriate representation of function changes.

The optimized designs (obtained for $p$ 3%) are collected in Table 3, whereas Fig. 3 shows the effects of selected FD steps on the TR performance for $x_1$ and $x_2$ designs. It should also be noted that the FD-induced change of computational cost is up to 22% for the considered tests, with the better performing setups being (understandably) more expensive, as more iterations are needed to accurately exploit the search space.

**Table 1.** Antenna designs used for TR algorithm benchmark.

| | Designs used for benchmark | | | | | | | | | |
|---|---|---|---|---|---|---|---|---|---|---|
| | $x_1^{(0)}$ | $x_2^{(0)}$ | $x_3^{(0)}$ | $x_4^{(0)}$ | $x_5^{(0)}$ | $x_6^{(0)}$ | $x_7^{(0)}$ | $x_8^{(0)}$ | $x_9^{(0)}$ | $x_{10}^{(0)}$ |
| $U(x)$ | –2.03 | –4.50 | –0.34 | –3.53 | –2.97 | –0.87 | –2.04 | –2.45 | –6.23 | –1.23 |
| Parameters | 17.5 | 22.2 | 18.8 | 17.9 | 14.6 | 13.4 | 20.4 | 13.6 | 18.2 | 21.6 |
| | 15.1 | 21.3 | 16.7 | 15.6 | 11.2 | 9.58 | 18.9 | 9.87 | 15.9 | 20.5 |
| | 6.79 | 8.79 | 7.30 | 6.95 | 5.52 | 4.99 | 8.04 | 5.08 | 7.07 | 8.56 |
| | 1.72 | 2.64 | 1.96 | 1.79 | 1.13 | 0.89 | 2.30 | 0.93 | 1.85 | 2.54 |
| | 9.07 | 12.8 | 10.0 | 9.37 | 6.73 | 5.75 | 11.3 | 5.92 | 9.60 | 12.3 |
| | 6.05 | 8.51 | 6.68 | 6.25 | 4.49 | 3.83 | 7.59 | 3.95 | 6.40 | 8.23 |

**Table 2.** TR-based optimization – benchmark results.

| Design | Perturbation size / type | | | | | | | | | | | | | | | |
|---|---|---|---|---|---|---|---|---|---|---|---|---|---|---|---|---|
| | $0.005x^{(0)}$ | | $0.01x^{(0)}$ | | $0.015x^{(0)}$ | | $0.02x^{(0)}$ | | $0.025x^{(0)}$ | | $0.03x^{(0)}$ | | $\varepsilon^{0.5}$ & | | Custom | |
| | $U(x^*)$ | $R^\#$ | $U(x^*)$ | $R^\#$ | $U(x^*)$ | $R^\#$ | $U(x^*)$ | $R^\#$ | $U(x^*)$ | $R^\#$ | $U(x^*)$ | $R^\#$ | $U(x^*)$ | $R^\#$ | $U(x^*)$ | $R^{\#,!}$ |
| $x_1$ | –16.9 | 74 | –19.8 | 67 | –17.3 | 45 | –23.8 | 88 | –22.1 | 88 | –26.9 | 88 | –21.2 | 88 | –26.7 | 104 |
| $x_2$ | –27.2 | 35 | –31.1 | 58 | –24.1 | 39 | –29.6 | 43 | –30.5 | 29 | –31.2 | 43 | –28.0 | 45 | –24.3 | 61 |
| $x_3$ | –29.7 | 72 | –30.2 | 86 | –30.6 | 88 | –30.2 | 86 | –25.3 | 45 | –29.5 | 81 | –27.9 | 77 | –24.3 | 92 |
| $x_4$ | –30.7 | 52 | –30.8 | 52 | –28.1 | 52 | –27.5 | 36 | –30.9 | 51 | –30.9 | 85 | –30.2 | 81 | –27.6 | 61 |
| $x_5$ | –16.3 | 60 | –17.0 | 30 | –18.4 | 44 | –17.6 | 39 | –21.5 | 94 | –28.2 | 94 | –17.7 | 50 | –15.4 | 40 |
| $x_6$ | –27.1 | 88 | –26.3 | 72 | –31.0 | 63 | –30.2 | 65 | –28.6 | 88 | –30.6 | 94 | –27.0 | 81 | –29.6 | 95 |
| $x_7$ | –20.2 | 88 | –11.4 | 59 | –26.7 | 66 | –29.5 | 93 | –28.7 | 94 | –29.2 | 94 | –28.7 | 94 | –20.4 | 63 |
| $x_8$ | –27.0 | 72 | –26.8 | 38 | –29.7 | 58 | –23.3 | 37 | –29.4 | 51 | –30.6 | 72 | –27.9 | 36 | –27.3 | 76 |
| $x_9$ | –26.6 | 35 | –27.8 | 42 | –28.1 | 66 | –29.6 | 65 | –30.1 | 72 | –26.3 | 37 | –17.9 | 37 | –30.0 | 68 |
| $x_{10}$ | –16.1 | 80 | –16.4 | 81 | –16.3 | 74 | –16.4 | 81 | –16.0 | 74 | –16.3 | 73 | –17.0 | 64 | –16.3 | 84 |
| $E^\$$ | –23.8 | 65 | –23.8 | 59 | –25.0 | 60 | –25.8 | 63 | –26.3 | 69 | –28.0 | 76 | –24.3 | 65 | –24.2 | 74 |
| $\sigma^\$$ | 5.76 | 20 | 7.03 | 18 | 5.68 | 15 | 5.27 | 23 | 4.96 | 23 | 4.43 | 20 | 5.25 | 22 | 5.21 | 19 |

$\$$  Average and standard deviation of the values obtained for the considered designs
$\#$  Computational cost of the optimization expressed in the number of EM model evaluations
$\&$  Perturbation size determined as square root of the machine precision $\varepsilon = 10^{-7}$
$!$  Including the cost of manual adjustment of perturbations

**Table 3.** Designs optimized using FD setup with 3 percent perturbations.

| | Optimized designs[$] – $p = 0.003x$ | | | | | | | | | |
|---|---|---|---|---|---|---|---|---|---|---|
| | $x_1^*$ | $x_2^*$ | $x_3^*$ | $x_4^*$ | $x_5^*$ | $x_6^*$ | $x_7^*$ | $x_1^*$ | $x_9^*$ | $x_{10}^*$ |
| Parameters | 20.2 | 20.7 | 20.5 | 20.6 | 20.4 | 20.5 | 20.5 | 20.6 | 19.9 | 19.0 |
| | 12.9 | 11.3 | 12.6 | 11.9 | 11.7 | 12.6 | 11.9 | 12.6 | 12.6 | 11.5 |
| | 4.10 | 3.56 | 3.92 | 3.72 | 3.84 | 3.84 | 3.73 | 3.85 | 4.32 | 5.22 |
| | 0.20 | 0.49 | 0.21 | 0.36 | 0.39 | 0.23 | 0.38 | 0.23 | 0.20 | 2.59 |
| | 10.1 | 11.1 | 10.4 | 10.8 | 10.4 | 10.8 | 11.0 | 10.7 | 9.70 | 3.00 |
| | 9.37 | 10.0 | 9.81 | 9.89 | 10.0 | 9.64 | 9.79 | 9.71 | 10.0 | 5.26 |

[$] For objective function values see Table 2

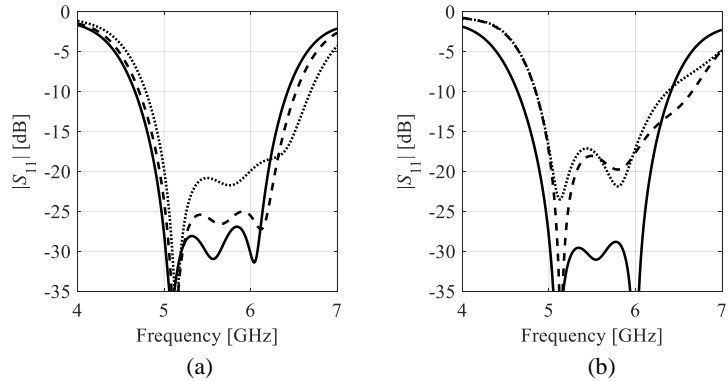

**Fig. 3.** Antenna design – comparison of the final designs obtained using TR-based optimization with forward FD steps of 1% (···), 2% (– –), and 3% (——) w.r.t. the center design: (a) design case $x_1$ and (b) design case $x_2$. Note that the selected FD perturbations have a notable effect on the quality of the optimized antenna response.

## 5    Conclusion

In this work, the effects of FD setup on the performance of TR-based antenna optimization have been demonstrated on a case study basis concerning a planar radiator. A total of 80 numerical experiments, spanned across 10 designs (each optimized using 8 different setups) have been performed. The numerical results demonstrate that the selected FD substantially affects the quality of the obtained designs (by as much as 18 dB per problem and around 4 dB on average), but also design cost (by up to 22%). Furthermore, the gathered—problem-specific data—suggest that determination of FD steps in a range between 1.5% and 3% of the initial design is suitable for ensuring high quality results for a range of test cases.

**Acknowledgments.** This work was supported in part by the National Science Centre of Poland Grants 2020/37/B/ST7/01448 and 2021/43/B/ST7/01856, National Centre for Research and Development Grant NOR/POLNOR/HAPADS/0049/2019-00, and Gdansk University of Technology (Excellence Initiative - Research University) Grant 16/2023/IDUB/IV.2/EUROPIUM.

**Disclosure of Interests.** The authors declare no conflicts of interest.